\documentclass[10pt]{article}

\usepackage{amsmath}
\usepackage{amssymb}
\usepackage{indentfirst}
\usepackage{graphics} 
\usepackage{color}

\makeatletter

\@addtoreset{equation}{section}
\makeatother
\def\<{\langle}
\def\>{\rangle}

\newtheorem{lem}{Lemma}[section]
\newtheorem{theo}{Theorem}[section]
\newtheorem{rem}{Remark}[section]

\makeatletter
   
   \@addtoreset{equation}{section}
\makeatother

\begin{document}
\title{\bf $L^{2}$-growth property for wave equations with\\ higher derivative terms}
\author{Ryo Ikehata\thanks{ikehatar@hiroshima-u.ac.jp} \\ {\small Department of Mathematics, Division of Educational Sciences}\\ {\small Graduate School of Humanities and Social Sciences} \\ {\small Hiroshima University} \\ {\small Higashi-Hiroshima 739-8524, Japan} \\and\\ Xiaoyan Li \thanks{Corresponding author: xiaoyanli@hust.edu.cn} \\ {\small School of Mathematics and Statistics}\\ {\small Huazhong University of Science and Technology} \\ {\small Wuhan, Hubei 430074, PR China}}
\date{}
\maketitle
\begin{abstract}
We consider the Cauchy problems in ${\bf R}^{n}$ for wave equations with higher derivative terms. We derive sharp growth estimates of the $L^{2}$-norm of the solution itself for the case of $n = 1$ and $n = 2$. By imposing the weighted $L^{1}$-initial velocity, we can get the lower and upper bound estimates of the solution itself.  For the case of $n\geq3$,  we observe that the $L^{2}$-growth behavior of the solution never occurs in the $(L^{2}\cap L^{1})$-framework of the initial data.
\end{abstract}
\section{Introduction}
\footnote[0]{Keywords and Phrases: Wave equation; higher derivative terms; weighted $L^{1}$-data; low dimensional case; growth estimates.}
\footnote[0]{2010 Mathematics Subject Classification. Primary 35L05; Secondary 35B40, 35C20, 35E15.}

We consider the Cauchy problem of the wave equation with a higher derivative term:
\begin{align}
& u_{tt} - \Delta u - \Delta u_{tt} = 0,\ \ \ (t,x)\in (0,\infty)\times {\bf R}^{n},\label{eqn}\\
& u(0,x)= u_0(x), \quad  u_{t}(0,x)= u_{1}(x),\ x\in{\bf R}^{n}.\label{initial}
\end{align}
Here, we assume, for the moment, $[u_{0},u_{1}] \in H^{1}({\bf R}^{n})\times H^{1}({\bf R}^{n})$.\\ 
\noindent
Concerning the existence of a unique energy solution to problem \eqref{eqn}-\eqref{initial}, by the Lumer-Phillips Theorem one can find that the problem (1.1)-(1.2) has a unique mild solution
\[u \in  C^{1}([0,\infty);H^{1}({\bf R}^{n})) \]
satisfying the energy conservation law such that
\begin{equation}\label{i-1}
E(t) = E(0),\quad t \geq 0,
\end{equation}
where the total energy $E(t)$ for the solution to problem \eqref{eqn}-\eqref{initial} can be defined by
\[E(t) := \frac{1}{2}\left(\Vert u_{t}(t,\cdot)\Vert_{L^{2}({\bf R}^{n})}^{2} + \Vert \nabla u_{t}(t,\cdot)\Vert_{L^{2}({\bf R}^{n})}^{2} + \Vert\nabla u(t,\cdot)\Vert_{L^{2}({\bf R}^{n})}^{2}\right).\] 
References \cite{H} and \cite{Mi-2} are helpful in detailing these discussions in Section 2 below.

The Strichartz and $L^p$-$L^q$ estimates for (free) wave equations
\[u_{tt}-\Delta u = 0\]
are powerful and particularly well-known as estimation formulas for the solution of the free wave equation itself (see \cite{S, W, B, P, M, Pl} and the references therein), however when one tries to obtain the $L^2$ estimate of the solution itself from some motivation, one will feel that it may be critical and a bit far from the best estimate, especially in low dimensions. Furthermore, it is important to recognize that the low-dimensional case, coupled with the infeasibility of the Hardy-type inequality, requires more delicate treatment than the high-dimensional case. Among such series of estimations, most of $L^2$ estimates, in some sense, are critical and difficult to derive. On the other hand, since the optimal growth in time estimate of the $L^2$ norm of the solution of the wave equation with strong damping in low dimensions was recently obtained in \cite{I-14, IO} to the equation 
\begin{equation}\label{ike-1000}
u_{tt}-\Delta u -\Delta u_{t} = 0,
\end{equation}
studies observing the best growth estimate of various wave-derived equations have been reported, one after another (see \cite{Ba, CI, CI-2, CT, FIM, IT, Mi}). Among them, Ikehata's estimates of the best $L^2$ norm for free wave and plate solutions are one of the most important fundamental results of these studies (see \cite{JHDE-ike, FE-ike}). The equation \eqref{ike-1000} expresses the wave part $+$ third-order derivative term in a sense. The purpose of this study is to contribute the $L^2$ estimates of the solution itself when higher derivative terms are added to the free wave equation. We believe that the results we claimed below are completely novel for this type equation. Incidentally, there has been some references \cite{DL, FIM-1}  related the $L^2$ behavior of the solution itself of the plate equation with higher derivative terms and dissipative terms. Although the form of the equation changes slightly, there are a lot of interesting literature recently on the subject of infinite-time blowup, and some of which are listed as \cite{CI, CL, CT-2, I-AA, ARR}.

By the way, although the equation (1.1) is related to the so-called generalized IMBq equation:
\begin{equation}\label{0-1}
u_{tt}-\Delta u_{tt}-\Delta u = \Delta f(u),
\end{equation}
it is natural to image that the result for the linear homogeneous case $f(u) = 0$ in particular has some influence on the behavior of the solution of the equation \eqref{0-1}, for example, a scattering result as in \cite[Theorem 1.2]{WC}, and thus it would be well worth considering (1.1). The relevant researches for equation \eqref{0-1} should be explored in \cite{WC-1, WC} and the references therein. Furthermore, for a study of the coupled Schr\"odinger equation and IMBq, we refer to  \cite{A}, \cite{OT} and the references therein. It would be interesting to explore whether or not there is an effect of singularity as found from (1.1) in this paper when coupling, but this is a future issue. It should be emphasized that as can be seen from some of the preceding papers (\cite{WC} and the references therein), it seems to make particular sense to deal with lower dimensions from a physics perspective as well.\\

 Before going to introduce our theorems, we present the following notations.
\\

{\bf Notation.} {\small Throughout this paper, $\| \cdot\|_q$ stands for the usual $L^q({\bf R}^{n})$-norm. For simplicity of notation, in particular, we use $\| \cdot\|$ instead of $\| \cdot\|_2$. 
	We also introduce the following weighted functional spaces
	\[L^{1,\gamma}({\bf R}^{n}) := \left\{f \in L^{1}({\bf R}^{n}) \; \bigm| \; \Vert f\Vert_{1,\gamma} := \int_{{\bf R}^{n}}(1+\vert x\vert^{\gamma})\vert f(x)\vert dx < +\infty\right\}.\]
The Fourier transform ${\cal F}_{x\to\xi}(f)(\xi)$ of $f(x)$ is defined by 
	\[{\cal F}_{x\to\xi}(f)(\xi) = \hat{f}(\xi) := \displaystyle{\int_{{\bf R}^{n}}}e^{-ix\cdot\xi}f(x)dx, \quad \xi \in {\bf R}^n,\]
	as usual with $i := \sqrt{-1}$, and ${\cal F}_{\xi\to x}^{-1}$ expresses its inverse Fourier transform. We denote the surface area of the $n$-dimensional unit ball by $\omega_{n} := \displaystyle{\int_{\vert\omega\vert = 1}}d\omega$.
For each $n = 1,2$, we set
\[I_{0,n} := \Vert u_{1}\Vert_{L^{2}({\bf R}^{n})} + \Vert u_{1}\Vert_{L^{1}({\bf R}^{n})}.\]
 }
For observing an essential part of growth property of the solution itself, we treat only the trivial initial amplitude case $u_{0}(x) \equiv 0$. 

Our first result is concerned with the optimal growth property in the case of $n = 1$.

\begin{theo}\label{theorem1}
Let $n = 1$, $u_{0} = 0$ and $u_{1} \in H^{1}({\bf R})$. Then, the solution $u(t,x)$ to problem \eqref{eqn}-\eqref{initial} satisfies the following properties under the additional regularity on the initial data:
\[ \Vert u(t,\cdot)\Vert_{L^{2}({\bf R})} \leq C_{1}I_{0,1}\sqrt{t},~~~~~\text{if}~~
u_{1} \in L^{1}({\bf R}),\]
\[ C_{2}\left\vert\int_{{\bf R}}u_{1}(x)dx\right\vert\sqrt{t} \leq \Vert u(t,\cdot)\Vert_{L^{2}({\bf R})}, ~~~~~\text{if}~~
u_{1} \in L^{1,\gamma}({\bf R}),~\gamma \in (\frac{1}{2},1],\]
for $t \gg 1$, where $C_{j} > 0$ {\rm ($j = 1,2$)} are constants depending only on the space dimension and $\gamma \in (\frac{1}{2},1]$.
\end{theo}
Our next result is the case of $n = 2$.
\begin{theo}\label{theorem2}
Let $n = 2$, and $u_{0} = 0$ and $u_{1} \in H^{1}({\bf R}^{2})$. Then, the solution $u(t,x)$ to problem \eqref{eqn}-\eqref{initial} satisfies the following properties under the additional regularity on the initial data:
\[ 
\Vert u(t,\cdot)\Vert_{L^{2}({\bf R}^{2})} \leq C_{1}I_{0,2}\sqrt{\log t}, ~~~\text{if}~~u_{1} \in L^{1}({\bf R}^{2}),\]
\[
 C_{2}\left\vert\int_{{\bf R}^{2}}u_{1}(x)dx\right\vert\sqrt{\log t} \leq \Vert u(t,\cdot)\Vert_{L^{2}({\bf R}^{2})},     ~~~~~\text{if}~~
 u_{1} \in L^{1,\gamma}({\bf R}^{2}),~\gamma \in (0,1],
  \]
for $t \gg 1$, where $C_{j} > 0$ {\rm ($j = 1,2$)} are constants depending only on the space dimension and $\gamma$. 
\end{theo}
\begin{rem}{\rm From these results above, it seems quite natural to choose initia data such as $(-\Delta)^{-\frac{1}{2}}u_{1} \in L^{2}({\bf R}^{n})$ in (for example) \cite{WC} to get global in time solutions of the equation \eqref{0-1}.}
\end{rem}
\begin{rem}{\rm The problem (1.1)-(1.2) is already studied in \cite[(2.23) of Theorem 2.1]{WC}, there the so-called $L^{2}$-$L^{2}$ bounded estimate of the solution to problem (1.1)-(1.2) is discussed. Our growth estimate just derived comes from $(L^{2}\cap L^{1})$-$L^{2}$ type estimate of the solution.}
\end{rem}

Let us explain where the unique difficulty of this problem arises compared to previous studies in observing the $L^2$ estimates of the solution itself. In particular, the difficulty is more pronounced in the two-dimensional treatment.
As in the usual treatment, in the Fourier space ${\bf R}_{\xi}^{n}$ the problem \eqref{eqn}-\eqref{initial} and its solution $u(t,x)$ can be transformed into the following ODE with parameter $\xi \in {\bf R}_{\xi}^{n}$
\begin{align}
& (1+\vert\xi\vert^{2})w_{tt} + |\xi|^2 w = 0,\ \ \ t>0,\quad \xi \in {\bf R}_{\xi}^{n},\label{5}\\
& w(0,\xi) = 0, \quad  w_{t}(0,\xi)= w_{1}(\xi),\ \ \ \xi \in{\bf R}^{n} ,\label{6}
\end{align}
where $w_{1}(\xi) := \hat{u}_1(\xi)$ and $w(t,\xi) := \hat{u}(t,\xi)$. Moreover, one can easily solve the problem \eqref{5}-\eqref{6} (formally) as follows:
\begin{equation}\label{7}
w(t,\xi) = \frac{\sin(tf(\vert\xi\vert))}{f(\vert\xi\vert)}w_{1}(\xi),
\end{equation}
and
\begin{equation}\label{8}
f(r) := \frac{r}{\sqrt{1+r^{2}}}.
\end{equation}
First, note that the proof of Theorem \ref{theorem1} in one dimension can be handled with the same strategy as the method in \cite{JHDE-ike}, but in the proof of Theorem 1.2 in two dimensions, the method in \cite{JHDE-ike} can not be applied directly. In \cite{JHDE-ike}, it uses the fact that the range of  function $f(r) = r$ is the half-space $[0,\infty)$ because there the free wave equation is studied. While, the range of $f(r)$ determined in \eqref{8} is in the bounded interval $[0,1)$, which is one of the factors requiring a decidedly different treatment and thus a new problem arises.

In conclusion, even with the addition of higher derivative terms, as in the free wave case, a certain singularity, expressed in augmented estimates, is included in the solution itself for lower dimensions $n = 1,2$. Therefore, it must be handled with sufficient delicacy. For the proof of Theorem \ref{theorem2}, we only use the method coming from \cite{I-04, I-14} and the integration by parts. The use of integration by parts is inspired by \cite[Proposition A.1.]{CT}, which is an improvement of \cite{JHDE-ike, FE-ike}.\\

The following three basic facts will be used throughout this paper.

We set (possibly $L = 1$)
\begin{equation}\label{i-2}
L := \sup_{\theta \ne 0}\left\vert \frac{\sin\theta}{\theta}\right\vert < +\infty.
\end{equation}
Furthermore, let $\delta_{0} \in (0,1)$ be a real number such that
\begin{equation}\label{i-3}
\left\vert \frac{\sin\theta}{\theta}\right\vert \geq \frac{1}{2}
\end{equation}
for all $\theta \in (0,\delta_{0}]$. We also prepare the fundamental inequality
\begin{equation}\label{i-4}
\vert a + b\vert^{2} \geq \frac{1}{2}\vert a\vert^{2} - \vert b\vert^{2} 
\end{equation}
for all $a, b \in {\bf C}$.\\

The paper is organized as follows. In Section 2, the well-posedness of equation \eqref{eqn}-\eqref{initial} will be showed. In Section 3,  we derive the lower bound estimates of the $L^{2}$-norm of solutions, and in Section 4 we obtain the upper bound estimates of the $L^{2}$-norm of solutions. By combining the results obtained in Sections 3 and 4, Theorems 1.1 and 1.2 are proved at a stroke. In Section 5, we consider the higher dimensional case as an additional remark.\\

\section{The well-posedness of the solution}
This section is concerned mainly with the well-posedness for the problem \eqref{eqn}-\eqref{initial}.  
A natural tendency to cope with this problem is to adopt semigroup theory. For the reader's convenience, we outline the proof based on the ideas coming from \cite{H} and \cite{Mi-2}.  

We denote $v=u_t$ and $A=-\Delta$. It follows from \eqref{eqn} that 
$$(I+A)v_t=-Au,$$
where $I$ is  identity operator in $H^1{({\bf R})^n}$.

Setting 
\begin{gather*}
U=	\begin{pmatrix}
		u\\v
	\end{pmatrix},~~~
U_0=	\begin{pmatrix}
	u_0\\u_1
\end{pmatrix},~~~
\mathcal{A} =	\begin{pmatrix}
	0& I \\ -P& 0
\end{pmatrix},~~~
P=-(I+A)^{-1}A, 
\end{gather*}
we have
\begin{equation} \label{}
	\left\{
	\begin{aligned}
		&\frac{d}{dt}U=\mathcal{A}U,\\
		&U(0,x) =U_0.
	\end{aligned}
	\right.
\end{equation}
Here $\mathcal{D}(P)$ is defined by 
\begin{align}
	\mathcal{D}(P)
	=\{&u\in H^1({\bf R}^n): \text{There exists} ~y_u\in H^1({\bf R}^n) 
	~\text{such that} ~ \notag \\  
	&~~~~~~~~~~~~~~~~~~(A^{\frac{1}{2}}u, A^{\frac{1}{2}}\phi )=(A^{\frac{1}{2}}y_u, A^{\frac{1}{2}}\phi )+(y_u, \phi),~~\forall \phi \in  H^1({\bf R}^n)
\} , \label{c1222}
\end{align}
where $(\cdot, \cdot )$ denotes the inner product of $L^2({\bf R}^n)$: 
$$(f,g)=\int_{{\bf R}^{n}} f(x)g(x) dx,~~~f,g\in L^2({\bf R}^n).$$ 
Note  that $\mathcal{D}(P)$ is not empty because $0\in \mathcal{D}(P)$ when we take $y=0$.  If $u\in \mathcal{D}(P)$, there exists a unique $y_u\in H^1({\bf R}^n)$ such that \eqref{c1222} holds. Otherwise, at least there exists  $y_u^1$ and $y_u^2 \in H^1({\bf R}^n)$ satisfying 
\begin{align}
	(A^{\frac{1}{2}}u, A^{\frac{1}{2}}\phi )=(A^{\frac{1}{2}}y_u^1, A^{\frac{1}{2}}\phi )+(y^1_u, \phi)  \label{c1223} \\ 
	(A^{\frac{1}{2}}u, A^{\frac{1}{2}}\phi )=(A^{\frac{1}{2}}y_u^2, A^{\frac{1}{2}}\phi )+(y^2_u, \phi) \label{c1224}
\end{align}
for every $\phi \in  H^1({\bf R}^n)$. Taking $z=y_u^1-y_u^2$, combining \eqref{c1223} and\eqref{c1224}  yields
$$
(A^{\frac{1}{2}}z, A^{\frac{1}{2}}\phi )+(z, \phi)=0.
$$
Let $\phi=z$, then we conclude that $z=0$ in $H^1({\bf R}^n)$.

The above arguments imply that the linear operator $P: u\rightarrow y_u $ is well defined for each $u \in \mathcal{D}(P)$.\\

Next we prove the fact $\mathcal{D}(P)=H^1({\bf R}^n)$, that is, for every $u\in H^1({\bf R}^n)$, there exists $y_u\in H^1({\bf R}^n)$ such that \eqref{c1222} holds.

First, for each $u\in H^1({\bf R}^n)$, we define the bounded and linear functional $F_u: H^1({\bf R}^n)\rightarrow {\bf R}$ as following 
\begin{equation}\label{c1225}
	<F_u, \phi>=(A^{\frac{1}{2}}u, A^{\frac{1}{2}}\phi),~~~~  \forall \phi \in H^1({\bf R}^n).
\end{equation}
By Riesz representation theorem, there exists a unique $y_u\in H^1({\bf R}^n)$ such that 
\begin{equation}\label{c1226}
<F_u, \phi>=(y_u, \phi)_{H^1({\bf R}^n)},~~~~  \forall \phi \in H^1({\bf R}^n),	
\end{equation}
where $(\cdot, \cdot)_{H^1({\bf R}^n)}$ denotes the inner product of Sobolev space $H^1({\bf R}^n)$ and it is equivalent to
\begin{equation}\label{c1227}
(f, ~g)_{H^1({\bf R}^n)}=(A^{\frac{1}{2}}f, A^{\frac{1}{2}}g )+(f, g).
\end{equation}
It follows from \eqref{c1225}, \eqref{c1226} and \eqref{c1227} that
\begin{equation}\label{1227}
	(A^{\frac{1}{2}}u, A^{\frac{1}{2}}\phi )=(A^{\frac{1}{2}}y_u, A^{\frac{1}{2}}\phi )+(y_u, \phi),~~\forall \phi \in  H^1({\bf R}^n).
\end{equation}
Therefore, we have $\mathcal{D}(P)=H^1({\bf R}^n)$ and we can define $P(u)=y_u$ on $H^1({\bf R}^n)$. Taking the Fourier transform of both sides of \eqref{1227} leads to 
\[
\int_{{\bf R}^n} 
|\xi|^2 \hat{u} \bar{\hat{\phi}} ~d\xi=\int_{{\bf R}^n} (|\xi|^2 +1)\hat{y}_u \bar{\hat{\phi}} ~ d\xi.
\]
Due to the arbitrariness of $\phi$, we obtain $|\xi|^2 \hat{u} =(|\xi|^2 +1)\hat{y}_u$, and then 

\begin{equation}\label{c1228}
	\hat{y}_u=\frac{|\xi|^2}{1+|\xi|^2} \hat{u}.
\end{equation}

Before giving the following lemma, we have to define the Hilbert space 
 $$\mathcal {H}:=H^1({\bf R}^n)\times H^1({\bf R}^n)$$
equipped with the inner product 
\[
\big( [y_1, z_1], [y_2, z_2]\big)_\mathcal{H} := (A^{\frac{1}{2}}y_1, A^{\frac{1}{2}}y_2 )+(y_1, y_2 )+(A^{\frac{1}{2}}z_1, A^{\frac{1}{2}}z_2 )+(z_1, z_2 ).
\]
\begin{lem} \label{c1235}
The operator 
\begin{gather*}
\mathcal{A} =	\begin{pmatrix}
	0& I \\ -P& 0
\end{pmatrix}: 
\mathcal {H} \rightarrow \mathcal{H}
\end{gather*}
generates a strongly continuous contraction semigroup $T(t)$ on $\mathcal {H}$.
\end{lem}
\noindent{ \it{ proof.}}
	For each $U=[u, v] \in H^1({\bf R}^n)\times H^1({\bf R}^n)$, it follows from \eqref{c1228} that
\begin{align}
    \big( \mathcal{A}U, U   \big)_\mathcal{H}&=\big( [v, -P(u)],~[u,v] \big)_\mathcal{H} \notag \\
	&=(A^{\frac{1}{2}}v, A^{\frac{1}{2}}u )+(v, u )-(A^{\frac{1}{2}}P(u), A^{\frac{1}{2}}v )-(P(u), v)\notag\\
	&=\int_{{\bf R}^n}(1+ |\xi|^2) \hat{v}\bar{\hat{u}}~d\xi-
	\int_{{\bf R}^n}\frac{ |\xi|^4+|\xi|^2}{|\xi|^2 +1} \hat{u}\bar{\hat{v}}~d\xi  \notag\\
	&=\int_{{\bf R}^n}\hat{v}\bar{\hat{u}}~d\xi + \int_{{\bf R}^n} |\xi|^2(\hat{v}\bar{\hat{u}}-\hat{u}\bar{\hat{v}})d\xi\notag\\ 
	&=\int_{{\bf R}^n}\hat{v}\bar{\hat{u}}~d\xi +2i \int_{{\bf R}^n} |\xi|^2 \text{Im}(\hat{u}\bar{\hat{v}})~d\xi.
\end{align}
This yields
\[{\text Re}\big( \mathcal{A}U, U   \big)_\mathcal{H} = {\text Re} \int_{{\bf R}^n} \hat{v}\bar{\hat{u}}~d\xi \leq {\text Re}\big( \frac{1}{2}U, U   \big)_\mathcal{H},\]
which implies
\begin{equation}\label{ike-1001}
{\text Re}\big( (\mathcal{A}-\frac{1}{2}\mathcal{I})U, U   \big)_\mathcal{H} \leq 0,
\end{equation}
where $\mathcal{I}$ is identity operator in $\mathcal{H}$.
\eqref{ike-1001} indicates that the operator $\mathcal{B} := \mathcal{A}-\frac{1}{2}\mathcal{I}$ is m-dissipative in Hilbert $\mathcal {H}$.\\

Next we prove that $\frac{1}{2}\mathcal{I}-\mathcal{B} = \mathcal{I}-\mathcal{A}$ is surjective, that is, for any fixed $[f,g]\in \mathcal{H}$, there exists $[u,v]\in H^1({\bf R}^n)\times H^1({\bf R}^n) $ such that 
\begin{gather} \label{c1229}
		(\mathcal{I}-\mathcal{A})	\begin{pmatrix}
		u\\v
	\end{pmatrix}
	=	\begin{pmatrix}
		f\\g
	\end{pmatrix}.
\end{gather}
By definition of $\mathcal{A}$, \eqref{c1229} is equivalent to
	\begin{align}
		&u-v=f,   \label{c1230}\\
		&Pu+v=g.   \label{c1231}
\end{align}
Submitting \eqref{c1230} into \eqref{c1231} leads to
\begin{equation*}
	(P+I)u=f+g.
\end{equation*}
By definition of $P$, one has
\begin{equation}\label{c1232}
(2A+I)u=(A+I)(f+g).
\end{equation}
Multiplying both sides of \eqref{c1232} by test function $h\in H^1({\bf R}^n)$ and integrating over  ${\bf R}^n$ yield
\begin{equation}\label{c1236}
	\Lambda(u, h)=F(h),
\end{equation}
where 
$$\Lambda(u, h)=2(A^{\frac{1}{2}}u, A^{\frac{1}{2}}h)+(u,h),~~~F(h) =(A^{\frac{1}{2}}(f+g), A^{\frac{1}{2}}h)+(f+g, h).$$
We see that $\Lambda(u, h)$ is continuous, coercive bilinear form on $H^1({\bf R}^n)$  and $F(h)$ is 
bounded linear functional related to  $f+g$. 
By the Lax-Milgram theorem, there exists a unique $u\in H^1({\bf R}^n)$ such that \eqref{c1236} holds. 
Then $v$ can be determined by \eqref{c1230}.

Finally, we conclude by the Lumer-Phillips theorem applied to the operator $\mathcal{A} = \frac{1}{2}\mathcal{I}+\mathcal{B}$ defined on $\mathcal{H}$ that the operator $\mathcal{A}$ generates  a strongly continuous contraction semigroup on Hilbert space $\mathcal{H}$, which completes the proof because the operator $\frac{1}{2}\mathcal{I}$ is bounded in $\mathcal{H}$, and $\mathcal{B}$ is m-dissipative in $\mathcal{H}$ (see \cite[Theorem 6.4]{G}). 
\hfill
$\Box$

By semigroup theory of linear operators and Lemma \ref{c1235}, the well-posedness of the equation \eqref{eqn}-\eqref{initial} can be obtained directly.

\begin{theo}\label{}
	Supposing initial data $[u_0, u_1]\in H^1({\bf R}^n)\times H^1({\bf R}^n) $, the equation \eqref{eqn}-\eqref{initial} admits a unique mild solution denoted by
	$$ [u, \partial_t u ]=T(t)[u_0, u_1]$$
	with regularity
	$$u \in  C^1([0,\infty); H^1({\bf R}^n)).$$
\end{theo}

\section{$L^{2}$-lower bound estimates of the solutions}

In this section, we derive the lower bound estimates in the case of $n = 1,2$. In particular, the estimates from below in the two-dimensional case require a delicate discussion to avoid the  inherent difficulties of the problem.\\

Now, let us give the lower bound estimates for $\Vert w(t,\cdot)\Vert$ in the case of $n = 1$, where $w(t,\xi)$ is defined in \eqref{7}.\\

At this first moment, one can assume that the initial velocity $u_{1} \in C_{0}^{\infty}({\bf R})$ by density.\\

Let $t > \delta_{0}$, $n \geq 1$ and denote a subset $L_{0}$ of ${\bf R}_{\xi}^{n}$ by
\begin{equation}\label{ike-10}
L_{0} := \{\xi \in {\bf R}_{\xi}^{n}\,:\,\vert\xi\vert \leq \frac{\delta_{0}}{\sqrt{t^{2}-\delta_{0}^{2}}}\}.
\end{equation}
Note that the function $f(r)$ defined in \eqref{8} is monotone increasing in $[0,\infty)$, $f(0) = 0$,\\
$\{f(r)\,\vert\,0 \leq r < \infty\} = [0,1)$, and 
\[\xi \in L_{0} \quad \Longleftrightarrow \quad tf(\vert\xi\vert) \in [0,\delta_{0}].\]
Therefore, using \eqref{i-3} one can proceed the estimate of an essential part of the solution $w(t,\xi)$ for $n=1$ as follows:
\begin{equation}\label{9}
I_{l}(t) := \int_{L_{0}}\frac{\sin^{2}(tf(\vert\xi\vert))}{f^{2}(\vert\xi\vert)}d\xi \geq \frac{t^{2}}{4}\int_{L_{0}}d\xi  = \frac{t^{2}}{2}\frac{\delta_{0}}{\sqrt{t^{2}-\delta_{0}^{2}}} \geq C_{0}t
\end{equation}
for $t \gg 1$, where $C_{0} > 0$ is an universal constant. 

Let us decompose the initial data $w_{1}(\xi)$ in the Fourier space
\begin{equation}\label{9-1}
w_{1}(\xi) = P + (A(\xi)-iB(\xi)),\quad \xi \in {\bf R}_{\xi}^{n},\quad (n \geq 1)
\end{equation}
where
\[P := \int_{{\bf R}^{n}}u_{1}(x)dx,\]
\[A(\xi) := \int_{{\bf R}^{n}}(\cos(x\xi)-1)u_{1}(x)dx, \quad B(\xi) := \int_{{\bf R}^{n}}\sin(x\xi)u_{1}(x)dx.\]
It is known (see \cite{I-04}) that with some constant $M > 0$ one has
\begin{equation}\label{i-8}
\vert A(\xi)-iB(\xi)\vert \leq M\vert\xi\vert^{\gamma}\Vert u_{1}\Vert_{1,\gamma},\quad \xi \in {\bf R}_{\xi}^{n},
\end{equation}
when $u_{1} \in L^{1,\gamma}({\bf R}^{n})$ and $\gamma \in (0,1]$. Then, it follows from \eqref{i-4} and \eqref{9-1} with $n = 1$ that
\begin{align} \label{i-9}
	J_{1}(t)&: = \int_{L_{0}}\frac{\sin^{2}(tf(\vert\xi\vert))}{f^2(\vert\xi\vert)}\vert w_{1}(\xi)\vert^{2}d\xi \notag \\
	&\geq \frac{P^{2}}{2}\int_{L_{0}}\frac{\sin^{2}(tf(\vert\xi\vert))}{f^2(\vert\xi\vert)}d\xi - \int_{L_{0}}\vert A(\xi)-iB(\xi)\vert^{2}\frac{\sin^{2}(tf(\vert\xi\vert))}{f^2(\vert\xi\vert)}d\xi \notag \\
	&=: \frac{P^{2}}{2}I_{l}(t)-R_{l}(t).
\end{align}
By \eqref{i-8}, $R_{l}(t)$ is estimated as follows for $\gamma \in (1/2,1]$
\begin{align}	\label{10}
	R_{l}(t) &\leq M^{2}\Vert u_{1}\Vert_{1,\gamma}^{2}\int_{L_{0}}\frac{r^{2\gamma}}{f(r)^{2}}d\xi\notag  \\
	&= M^{2}\Vert u_{1}\Vert_{1,\gamma}^{2}\int_{0}^{\frac{\delta_{0}}{\sqrt{t^{2}-\delta_{0}^{2}}}}(1+r^{2})r^{2(\gamma-1)}dr\notag \\
		&\leq \frac{CM^{2}}{2\gamma-1}\Vert u_{1}\Vert_{1,\gamma}^{2}t^{-(2\gamma-1)},
\end{align}
where $r=|\xi|$, the constant $C > 0$, and $t \gg 1$. Therefore, from \eqref{9}, \eqref{i-9} and \eqref{10}, one can get the desired lower bound estimate
\[\Vert w(t,\cdot)\Vert^{2} \geq J_{1}(t) \geq C_{1}P^{2}t - C_{\gamma}\Vert u_{1}\Vert_{1,\gamma}^{2}t^{-(2\gamma-1)}, \quad t \gg 1\]
with some constants $C_{1} > 0$ and $C_{\gamma} > 0$. By density argument, one can state the following lemma.
\begin{lem}\label{lem1}Let $n = 1$, and $\gamma \in (\frac{1}{2},1]$. Assume $u_{1} \in L^{1,\gamma}({\bf R})$. Then, it holds that
\[\Vert w(t,\cdot)\Vert^{2} \geq CP^{2} t, \quad t \gg 1.\]
\end{lem}


Next, we discuss the proof of the two dimensional case. We rely on the integration by parts together with a trick initiated from \cite{JHDE-ike}.\\

Using \eqref{i-4} and \eqref{9-1}, we start with the integral combing the trick function $e^{-\vert\xi\vert^{2}}$

\begin{align}
	\Vert w(t,\cdot)\Vert^{2} &= \int_{{\bf R}^{2}}\frac{\sin^{2}(tf(\vert\xi\vert))}{f^{2}(\vert\xi\vert)}\vert w_{1}(\xi)\vert^{2}d\xi\notag\\
	&\geq \int_{{\bf R}^{2}}e^{-\vert\xi\vert^{2}}\frac{\sin^{2}(tf(\vert\xi\vert))}{f^{2}(\vert\xi\vert)}\vert P+(A(\xi)-iB(\xi))\vert^{2}d\xi\notag\\
	&\geq \frac{1}{2}P^{2}\int_{{\bf R}^{2}}e^{-\vert\xi\vert^{2}}\frac{\sin^{2}(tf(\vert\xi\vert))}{f^{2}(\vert\xi\vert)}d\xi - \int_{{\bf R}^{2}}e^{-\vert\xi\vert^{2}}\frac{\sin^{2}(tf(\vert\xi\vert))}{f^{2}(\vert\xi\vert)}\left(M^{2}\Vert u_{1}\Vert_{1,\gamma}^{2}\vert\xi\vert^{2\gamma}\right)d\xi\notag\\
	&\geq \frac{1}{2}P^{2}\int_{{\bf R}^{2}}e^{-\vert\xi\vert^{2}}\frac{\sin^{2}(tf(\vert\xi\vert))}{f^{2}(\vert\xi\vert)}d\xi - M^{2}\Vert u_{1}\Vert_{1,\gamma}^{2}\int_{{\bf R}^{2}}e^{-\vert\xi\vert^{2}}\frac{\vert\xi\vert^{2\gamma}}{f^{2}(\vert\xi\vert)}d\xi\notag\\
	&=:  \frac{1}{2}P^{2}T(t) -M^{2}\Vert u_{1}\Vert_{1,\gamma}^{2}U(t).\label{ike-51}
\end{align}

Incidentally, this crucial idea to use the trick function has been initiated in the paper \cite{JHDE-ike}. Of course, there are other ways to choose this function.

Now, $U(t)$ can be estimated as follows:
\begin{align*}
	\frac{U(t)}{\omega_{2}} &= \int_{0}^{\infty}e^{-r^{2}}r^{2\gamma-1}dr +  \int_{0}^{\infty}e^{-r^{2}}r^{2\gamma+1}dr\\
	&= \frac{1}{\gamma}\int_{0}^{\infty}e^{-r^{2}}r^{2\gamma+1}dr +  \int_{0}^{\infty}e^{-r^{2}}r^{2\gamma+1}dr\\
	&=: K_{0}.
\end{align*}
Note that $K_{0} > 0$ has a finite value for $\gamma \in (0,1]$. The role of trick function $e^{-\vert\xi\vert^{2}}$ is useful. Thus, from \eqref{ike-51} one has
\begin{equation}\label{ike-52}
\Vert w(t,\cdot)\Vert^{2} \geq \frac{1}{2}P^{2}T(t) -M^{2}\Vert u_{1}\Vert_{1,\gamma}^{2}\omega_{2}K_{0},\quad t > 0.
\end{equation}

Finally, let us estimate the main term $T(t)$.\\

Since
\[2\sin^{2}(tf(r)) = 1-\cos(2tf(r)),\]
we get for $t > 1$
\begin{align}
	T(t) &\geq \frac{\omega_{2}}{2}\int_{1/t}^{1}e^{-r^{2}}\frac{2 \sin^{2}(tf(r))}{f^{2}(r)} rdr\notag\\
	&= \frac{\omega_{2}}{2}\int_{1/t}^{1}e^{-r^{2}}\frac{r}{f^{2}(r)}dr - \frac{\omega_{2}}{2}\int_{1/t}^{1}e^{-r^{2}}\frac{r}{f^2(r)}\cos(2tf(r))dr\notag\\
	&=: \frac{\omega_{2}}{2}T_{1}(t) - \frac{\omega_{2}}{2}T_{2}(t).\label{ike-53}
\end{align}
Here, we see that
\begin{equation}\label{ike-54}
	T_{1}(t) = \int_{1/t}^{1}e^{-r^{2}}\frac{1+r^{2}}{r}dr \geq \int_{1/t}^{1}e^{-r^{2}}\frac{1}{r}dr
\geq e^{-1}\int_{1/t}^{1}\frac{1}{r}dr = e^{-1}\log t,\quad t > 1.
\end{equation}
\noindent
$T_{2}(t)$ can be decomposed into two parts:
\begin{align}\label{ike-55}
T_{2}(t) &= \int_{1/t}^{1}e^{-r^{2}}r^{-1}\cos(2tf(r))dr + \int_{1/t}^{1}e^{-r^{2}}r\cos(2tf(r))dr \notag \\
&=: R_{1}(t) + R_{2}(t).
\end{align}
$R_{2}(t)$ can be estimated easily by
\begin{equation}\label{ike-56}
\vert R_{2}(t)\vert \leq \int_{1/t}^{1}e^{-r^{2}}rdr = \frac{1}{2}(e^{-\frac{1}{t^{2}}}-e^{-1}) \leq \frac{1}{2}, \quad t > 1.
\end{equation}

In order to estimate $R_{1}(t)$, we develop the integration by parts. This idea is inspired by \cite[Proposition A.1]{CT}. Indeed, since 
\[\cos(2tf(r)) = \frac{1}{2f'(r)t}  \left( \frac{d}{dr}\sin(2tf(r))\right),\]
\[ f'(r) = \frac{1}{(1+r^{2})\sqrt{1+r^{2}}},\]
\begin{equation}\label{ike-57}
	R_{1}(t) = \frac{1}{2t}\int_{1/t}^{1}\frac{e^{-r^{2}}}{r}(1+r^{2})\sqrt{1+r^{2}}\frac{d}{dr}\sin(2tf(r))dr
=: \frac{1}{2t}K(t),
\end{equation}
$K(t)$ can be estimated by
\begin{align}
	K(t) &= \left[\frac{e^{-r^{2}}}{r}(1+r^{2})\sqrt{1+r^{2}}\sin(2tf(r))\right]_{1/t}^{1}\notag \\
	&~~~~-\int_{1/t}^{1}\frac{d}{dr}\left(\frac{e^{-r^{2}}}{r}(1+r^{2})^{\frac{3}{2}}\right)\sin(2tf(r))dr \notag\\
	&=: K_{1}(t) + K_{2}(t).\label{ike-58}
\end{align}
It is obvious that
\begin{equation}\label{ike-59}
\vert K_{1}(t)\vert \leq C_{1} + C_{2}t
\end{equation}
with some constants $C_{1} > 0$ and $C_{2} > 0$. Next, we should check the upper bound estimate for $\vert K_{2}(t)\vert$. In fact, it holds that
\[\frac{d}{dr}\left(\frac{e^{-r^{2}}}{r}(1+r^{2})^{\frac{3}{2}}\right) = -\frac{e^{-r^{2}}}{r^{2}}(2r^{2}+1)(1+r^{2})^{\frac{3}{2}} + 3e^{-r^{2}}(1+r^{2})^{\frac{1}{2}},\]
\begin{align}
K_{2}(t) &=	\int_{1/t}^{1}\frac{e^{-r^{2}}}{r^{2}}(2r^{2}+1)(1+r^{2})^{\frac{3}{2}}\sin(2tf(r))dr-3\int_{1/t}^{1}e^{-r^{2}}(1+r^{2})^{\frac{1}{2}}\sin(2tf(r))dr\notag\\
&=: K_{2,1}(t) -K_{2,2}(t).\label{ike-60}
\end{align}
Simple calculations yield
\begin{equation}\label{ike-61}
\vert K_{2,2}(t)\vert \leq 3\int_{1/t}^{1}e^{-r^{2}}(1+r^{2})^{\frac{1}{2}}dr \leq 3e^{-\frac{1}{t^{2}}}\sqrt{2}(1-\frac{1}{t}) \leq C_{3}
\end{equation}
for some constant $C_{3} > 0$. While,
\begin{equation}\label{ike-62}
\vert K_{2,1}(t)\vert \leq \int_{1/t}^{1}e^{-r^{2}}r^{-2}(2r^{2}+1)(1+r^{2})^{\frac{3}{2}}dr \leq 6\sqrt{2}e^{-\frac{1}{t^{2}}}\int_{1/t}^{1}r^{-2}dr \leq C_{4}t 
\end{equation}
for $t \gg 1$ and some constant $C_{4} > 0$.

By summerizing \eqref{ike-55}-\eqref{ike-62},  the estimate for $T_{2}(t)$ is given by
\begin{equation}\label{ike-63}
\vert T_{2}(t)\vert \leq \frac{1}{2t}(C_{1} + C_{2}t + C_{3} + C_{4}t), \quad t \gg 1. 
\end{equation}
Therefore, by \eqref{ike-52}, \eqref{ike-53}, \eqref{ike-54} and \eqref{ike-63}, one can state the following lemma.
\begin{lem}\label{lem2}Let $n = 2$, and $\gamma \in (0,1]$. Assume $u_{1} \in L^{1,\gamma}({\bf R}^{2})$. Then, it holds that
\[\Vert w(t,\cdot)\Vert^{2} \geq CP^{2}\log t, \quad t \gg 1.\]
\end{lem}


\section{$L^{2}$-upper bound estimates of the solution}

In this section, we derive upper bound estimate of $\Vert u(t,\cdot)\Vert$ as $t \to \infty$. We use the $L^{2}$-regularity together with $L^{1}$-regularity of the initial velocity.\\ 
We first treat the one dimensional case. It also suffices to derive the desired upper bound estimate by assuming $u_{1} \in C_{0}^{\infty}({\bf R}^n)$.\\

From \eqref{7} and \eqref{8}, we can proceed the estimate as follows:
\begin{align}
	\Vert w(t,\cdot)\Vert^{2} &= \int_{{\bf R}_{\xi}}\vert\frac{\sin(tf(\vert\xi\vert))}{f(\vert\xi\vert)}\vert^{2}\vert w_{1}(\xi)\vert^{2}d\xi\notag\\
	&= \int_{L_{0}}\frac{\sin^{2}(tf(\vert\xi\vert))}{f^{2}(\vert\xi\vert)}\vert w_{1}(\xi)\vert^{2}d\xi + \int_{{\bf R}_{\xi}\setminus L_{0}}\frac{\sin^{2}(tf(\vert\xi\vert))}{f^{2}(\vert\xi\vert)}\vert w_{1}(\xi)\vert^{2}d\xi\notag \\
	&=: L_{1}(t) + L_{2}(t), \label{i-14}
\end{align}
where the set $L_{0}$ with $n = 1$ is defined in \eqref{ike-10}. Then, from \eqref{i-2} it holds that
\begin{equation}\label{15}
L_{1}(t) \leq L^{2}t^{2}\omega_{1}\Vert u_{1}\Vert_{1}^{2}\int_{0}^{A(t)}dr = L^{2}t^{2}\omega_{1}\Vert u_{1}\Vert_{1}^{2}A(t) \leq C_{1}\Vert u_{1}\Vert_{1}^{2}t
\end{equation}
for $t \gg 1$, where $C_{1} > 0$ is a constant, and 
\begin{equation}\label{155-1}
A(t) := \frac{\delta_{0}}{\sqrt{t^{2}-\delta_{0}^{2}}}.
\end{equation}

While, for $L_{2}(t)$ one has
\begin{align}
L_{2}(t) &\leq \int_{{\bf R}_{\xi}\setminus L_{0}}\frac{1+\vert\xi\vert^{2}}{\vert\xi\vert^{2}}\vert w_{1}(\xi)\vert^{2}d\xi\notag\\
&= \int_{{\bf R}_{\xi}\setminus L_{0}}\frac{1}{\vert\xi\vert^{2}}\vert w_{1}(\xi)\vert^{2}d\xi + \int_{{\bf R}_{\xi}\setminus L_{0}}\vert w_{1}(\xi)\vert^{2}d\xi\notag\\
&\leq \omega_{1}\Vert u_{1}\Vert_{1}^{2}\int_{A(t)}^{\infty}r^{-2}dr + \Vert u_{1}\Vert^{2} = \omega_{1}\Vert u_{1}\Vert_{1}^{2}\frac{1}{A(t)} + \Vert u_{1}\Vert^{2}\notag\\
&\leq C_{2}\Vert u_{1}\Vert_{1}^{2}t + \Vert u_{1}\Vert^{2} \label{16}	
\end{align}
for $t \gg 1$, where $C_{2} > 0$ is a constant. Thus, from \eqref{15} and \eqref{16} one can derive the desired estimate as the following lemma.
\begin{lem}\label{lem3}Let $n = 1$, and $u_{1} \in L^{1}({\bf R})\cap L^{2}({\bf R})$. Then, it holds that
\[\Vert w(t,\cdot)\Vert^{2} \leq C(\Vert u_{1}\Vert^{2} + \Vert u_{1}\Vert_{1}^{2}t), \quad t \gg 1.\]
\end{lem}


Next, let us give the upper bound estimate for $n = 2$ based on the decomposition below
\begin{align}
\Vert w(t,\cdot)\Vert^{2} &= \int_{{\bf R}_{\xi}^{2}}\vert\frac{\sin(tf(\vert\xi\vert))}{f(\vert\xi\vert)}\vert^{2}\vert w_{1}(\xi)\vert^{2}d\xi \notag\\
&= \int_{L_{0}}\frac{\sin^{2}(tf(\vert\xi\vert))}{f^{2}(\vert\xi\vert)}\vert w_{1}(\xi)\vert^{2}d\xi + \int_{{\bf R}_{\xi}^{2}\setminus L_{0}}\frac{\sin^{2}(tf(\vert\xi\vert))}{f^{2}(\vert\xi\vert)}\vert w_{1}(\xi)\vert^{2}d\xi\notag\\
&=: G_{1}(t) + G_{2}(t), \label{17}	
\end{align}
where the set $L_{0}$ with $n = 2$ is defined in \eqref{ike-10}. It also suffices to assume $u_{1} \in C_{0}^{\infty}({\bf R}^{2})$ in the derivation of upper bound estimate.

We first treat $G_{1}(t)$ to get the estimate such that
\begin{equation}\label{18}
G_{1}(t) \leq L^{2}t^{2}\omega_{2}\Vert u_{1}\Vert_{1}^{2}\int_{0}^{A(t)}rdr = 2^{-1}L^{2}\omega_{2}\Vert u_{1}\Vert_{1}^{2}\delta_{0}^{2}\frac{t^{2}}{t^{2}-\delta_{0}^{2}} \leq C_{3}\Vert u_{1}\Vert_{1}^{2}
\end{equation}
with some constant $C_{3} > 0$ and $t \gg 1$, where one uses the fact  \eqref{i-2}.

For the estimate of $G_{2}(t)$, it must be decomposed into three parts as follows in order to get the desired growth estimate
\begin{equation}\label{19}
G_{2}(t) = g_{1}(t) + g_{2}(t) + g_{3}(t) := \left(\int_{\vert\xi\vert \geq C(t)} + \int_{C(t) \geq \vert\xi\vert \geq B(t)} + \int_{B(t) \geq \vert\xi\vert \geq A(t)}\right)\frac{\sin^{2}(tf(\vert\xi\vert))}{f^{2}(\vert\xi\vert)}\vert w_{1}(\xi)\vert^{2}d\xi,
\end{equation}
where $t \gg 1$, $A(t)$ is already defined in \eqref{155-1}, and 
\[B(t) := \frac{\delta_{0}}{\sqrt{t -\delta_{0}^{2}}},\]
\[C(t) := \frac{\delta_{0}}{\sqrt{\log t -\delta_{0}^{2}}}.\]

Now, we see that
\begin{equation}\label{20}
g_{1}(t) \leq \frac{1}{f^{2}(C(t))}\int_{ |\xi| \geq C(t)}\vert w_{1}(\xi)\vert^{2}d\xi \leq \frac{\log t}{\delta_{0}^{2}}\Vert u_{1}\Vert^{2}, \quad t \gg 1.
\end{equation}
Recall that the function $f(r)$ is monotone increasing on $[0,\infty)$.

Next, about $g_{2}(t)$ we have
\begin{align}
	g_{2}(t) &\leq \omega_{2}\Vert u_{1}\Vert_{1}^{2}\int_{B(t)}^{C(t)}\frac{1+r^{2}}{r}dr = \omega_{2}\Vert u_{1}\Vert_{1}^{2}\left[\log r + \frac{1}{2}r^{2}\right]_{r= B(t)}^{r=C(t)}\notag\\
	&= \omega_{2}\Vert u_{1}\Vert_{1}^{2}\left(\log C(t)-\log B(t) + \frac{1}{2}(C(t)^{2}-B(t)^{2})\right)\notag\\
	&= \omega_{2}\Vert u_{1}\Vert_{1}^{2}\left(\frac{1}{2}\log(t-\delta_{0}^{2}) - \frac{1}{2}\log(\log t-\delta_{0}^{2}) + \frac{\delta_{0}^{2}}{2}(\frac{1}{\log t - \delta_{0}^{2}}- \frac{1}{t - \delta_{0}^{2}})\right)\notag\\
	&\leq C_{4}\Vert u_{1}\Vert_{1}^{2}\log t \label{21}
\end{align}
	with some constant $C_{4} > 0$ and $t \gg 1$.

Thirdly, $g_{3}(t)$ can be estimated similarly to \eqref{21}
\begin{align}
g_{3}(t) &\leq \omega_{2}\Vert u_{1}\Vert_{1}^{2}\int_{A(t)}^{B(t)}\frac{1+r^{2}}{r}dr \notag\\
&= \omega_{2}\Vert u_{1}\Vert_{1}^{2}]\left[\log r + \frac{1}{2}r^{2}\right]_{r= A(t)}^{r=B(t)}\notag\\
&= \omega_{2}\Vert u_{1}\Vert_{1}^{2}\left(\frac{1}{2}\log(t^{2}-\delta_{0}^{2}) - \frac{1}{2}\log(t-\delta_{0}^{2}) + \frac{\delta_{0}^{2}}{2}(\frac{1}{t - \delta_{0}^{2}}- \frac{1}{t^{2} - \delta_{0}^{2}})\right)\notag\\
&\leq C_{5}\Vert u_{1}\Vert_{1}^{2}\log t \label{22}	
\end{align}
with some constant $C_{5} > 0$ and $t \gg 1$. 

Therefore, the following lemma is a direct consequence of  \eqref{17}, \eqref{18}, \eqref{19}, \eqref{20}, \eqref{21} and \eqref{22}.
\begin{lem}\label{lem44}Let $n = 2$, and $u_{1} \in L^{1}({\bf R}^{2})\cap L^{2}({\bf R}^{2})$. Then, it holds that
\[\Vert w(t,\cdot)\Vert^{2} \leq C(\Vert u_{1}\Vert^{2} +\Vert u_{1}\Vert_{1}^{2})\log t, \quad t \gg 1.\]
\end{lem}

\par
\vspace{0.5cm}
Finally, the proofs of Theorems \ref{theorem1} and \ref{theorem2} are direct consequences of Lemmas \ref{lem1}, \ref{lem2}, \ref{lem3} and \ref{lem44}, and the Plancherel Theorem.\\

\section{$L^{2}$-upper bound estimates for the case $n \geq 3$}

In this section, we give some comments in the case of $n \geq 3$. In this case, the $L^2$-increasing property never holds, which is verified in below by the following calculations so far.\\

In fact, from \eqref{7} and \eqref{8} we can proceed the estimate as follows:
\begin{align}
	\Vert w(t,\cdot)\Vert^{2} &= \int_{{\bf R}_{\xi}^{n}}\vert\frac{\sin(tf(\vert\xi\vert))}{f(\vert\xi\vert)}\vert ^{2}\vert w_{1}(\xi)\vert^{2}d\xi\notag\\
	&= \int_{M_{0}}\frac{\sin^{2}(tf(\vert\xi\vert))}{f^{2}(\vert\xi\vert)}\vert w_{1}(\xi)\vert^2 d\xi + \int_{{\bf R}_{\xi}^{n}\setminus M_{0}}\frac{\sin^{2}(tf(\vert\xi\vert))}{f^{2}(\vert\xi\vert)}\vert w_{1}(\xi)\vert^{2}d\xi\notag\\
	&=: M_{1}(t) + M_{2}(t), \label{i-14-1}
\end{align}
where the set $M_{0}$ with $n \geq 3$ is defined by
\begin{equation}
M_{0} := \{\xi \in {\bf R}_{\xi}^{n}\,:\,\vert\xi\vert \leq \frac{1}{\sqrt{3}}\}.
\end{equation}
Note that  
\[\xi \in M_{0} \quad \Longleftrightarrow \quad f(\vert\xi\vert) \leq \frac{1}{2}.\]
Then, it holds that
\begin{equation}\label{15-1}
M_{1}(t) \leq \int_{M_0}\frac{1}{f^{2}(\vert\xi\vert)}\vert w_{1}(\xi)\vert^{2}d\xi \leq \omega_{n}\Vert u_{1}\Vert_{1}^{2}\int_{0}^{\frac{1}{\sqrt{3}}}(1+r^{2})r^{n-3}dr  \leq C_{1,n}\Vert u_{1}\Vert_{1}^{2}
\end{equation}
for $t \geq 0$, where $C_{1,n} > 0$ is a constant. 

While, $M_{2}(t)$ can be estimated by
\begin{equation}\label{16-1}
M_{2}(t) \leq \int_{{\bf R}_{\xi}^{n}\setminus M_{0}}\frac{1}{f^{2}(\vert\xi\vert)}\vert w_{1}(\xi)\vert^{2}d\xi \leq 4\int_{{\bf R}_{\xi}^{n}\setminus M_{0}}\vert w_{1}(\xi)\vert^{2}d\xi\leq C_{2,n}\Vert u_{1}\Vert
\end{equation}
for $t \geq 0$, where $C_{2,n} > 0$ is a constant. Thus, by \eqref{i-14-1}, \eqref{15-1} and \eqref{16-1}, one can obtain the following desired estimate.
\begin{lem}\label{lem4}Let $n \geq 3$, and $u_{1} \in L^{1}({\bf R})\cap L^{2}({\bf R})$. Then, it holds that
\[\Vert u(t,\cdot)\Vert^{2} = \Vert w(t,\cdot)\Vert^{2} \leq C(\Vert u_{1}\Vert^{2} + \Vert u_{1}\Vert_{1}^{2}), \quad t \gg 1.\]
\end{lem}
\begin{rem}{\rm Even if non-trivial initial value $u_{0}(x)$ is added, the result remains the same because of the easy estimate:
\[\int_{{\bf R}^{n}}\cos^{2}(f(\xi)t)\vert\hat{u}_{0}(\xi)\vert^{2}d\xi \leq \int_{{\bf R}^{n}}\vert\hat{u}_{0}(\xi)\vert^{2}d\xi = \Vert u_{0}\Vert^{2}.\]
So, the growth property never occurs for $n \geq 3$. In this sense, $n = 1,2$ are exceptional.}
\end{rem}

\vspace{0.5cm}
\noindent{\em Acknowledgement.}
\smallskip
This paper was written during Xiaoyan Li's stay as an overseas researcher at Hiroshima University from 12 December, 2022 to 11 December, 2023 under Ikehata's supervision as a host researcher. The work of the first author (Ryo Ikehata) was supported in part by Grant-in-Aid for Scientific Research (C) 20K03682 of JSPS. This work of the second author (Xiaoyan Li) was financially supported in part by Chinese Scholarship Council (Grant No. 202206160071). 

\vspace{0.5cm}
\noindent{\bf Declarations}
\smallskip

\noindent{\bf Data availability}
\smallskip
Data sharing not applicable to this article as no datasets were generated or analysed during the current study.

\noindent{\bf Conflict of interest}
\smallskip
The authors declare that they have no conflict of interest.


\end{document}